# Circumferential Crack Modeling of Thin Cylindrical Shells in Modal Deformation

Ali Alijani[a,*], Olga Barrera[b], Stéphane P.A. Bordas[c,†]

[a] *Department of Mechanical Engineering, Bandar Anzali Branch, Islamic Azad University, Bandar Anzali, Iran. (ORCID: 0000-0001-7782-9026)*
[b] *School of Engineering, Computing and Mathematics, Oxford Brookes University, UK. (ORCID: 0000-0002-0077-9582)*
[c] *Institute of Computational Engineering, Faculty of Sciences and Technology, University of Luxembourg, Luxembourg City, Luxembourg. (ORCID: 0000-0001-8634-7002)*

## Abstract

An innovative technique, called *conversion*, is introduced to model circumferential cracks in thin cylindrical shells. The semi-analytical finite element method is applied to investigate the modal deformation of the cylinder. An element including the crack is divided into three sub-elements with four nodes in which the stiffness matrix is enriched. The crack characteristics are included in the finite element method relations through conversion matrices and a rotational spring corresponding to the crack. Conversion matrices obtained by applying continuity conditions at the crack tip are used to transform displacements of the middle nodes to those of the main nodes. Moreover, another technique, called *spring set*, is represented based on a set of springs to model the crack as a separated element. Components of the stiffness matrix related to the separated element are incorporated while the geometric boundary conditions at the crack tip are satisfied. The effects of the circumferential mode number, the crack depth and the length of the cylinder on the critical buckling load are investigated. Experimental tests, ABAQUS modeling and results from literature are used to verify and validate the results and derived relations. In addition, the crack effect on the natural frequency is examined using the vibration analysis based on the conversion technique.

## 1. Introduction

A number of numerical methods, including the finite element method (FEM), the extended finite element method (XFEM), Meshfree, etc., have been developed to analyze discontinuous structures, e.g. [1-7]. Discrete spring models are an alternative created based on relations between the energy release rate and the stress intensity factors [8-12]. This discrete spring model has been used to analyze different engineering problems of cracked beams see e.g. [13-

[*] alijani@iaubanz.ac.ir (A. Alijani).
  nimalijanimech@yahoo.com (A. Alijani).
[†] stephane.bordas@alum.northwestern.edu (S.P.A. Bordas)
Fax: +98 134 440 0486. Postal Code: 43131-11111.

15]. Alijani et al. [16-18] presented a novel technique in the finite element method to include cracks into a beam element. They enriched components of the stiffness matrices by using crack properties modeled as a rotational spring.

The modal deformation of shells is important in engineering applications. The semi-analytical finite element is an efficient method in modal deformation analysis. Alijani et al. [19] and [20] introduced a new semi-analytical nonlinear finite element formulation according to a continuum-based approach to analyze the post-buckling of thin cylindrical shells under mechanical and thermal loads. Akrami and Erfani [21] investigated the critical buckling load for a circumferentially cracked cylindrical shell in which the cylinder and the crack are modeled as a beam-column on an elastic foundation and the rotational spring, respectively. Delale and Erdogan [22] proposed an approximate solution for a cylindrical shell containing a part-through surface crack assumed as circumferential or semi-elliptic. Ezzat and Erdogan [23] compared the experimental and theoretical results after discussing the analytical techniques used in modeling the problem of fatigue crack propagation of a cylindrical shell containing a circumferential flaw. Moradi and Tavaf [24] used the differential quadrature method combined with an evolutionary optimization algorithm to detect the crack position in cylindrical shell structures, where a circumferential crack is modeled by a rotational spring. Naschie [25] represented an initial post-buckling analysis for a simply supported concrete cylinder containing a circumferential crack. An eigenvalue buckling analysis [26] was carried out to investigate the effects of various parameters of cracked functionally graded cylindrical shells in the framework of the extended finite element method. Natarajan et al. [27] and [28] represented numerical solution and advanced discretization techniques in the buckling analysis of discontinuous thin-walled structures. In the aforementioned research works, some methods including analytical, approximate and numerical have been used to address the buckling problem of discontinuous structures.

In the present paper, a semi-analytical finite element method is initially applied to involve a circumferential crack in a thin cylindrical shell. Two techniques are used to enrich the stiffness matrix. The first technique, which was already applied for the analysis of cracked beams called the *conversion* technique [16-18], has been originally implemented to formulate the finite element relations for the cracked cylindrical shell. The second technique, called *spring set*, is applied by considering the crack as an element, whose stiffness matrix is assembled with the standard stiffness matrices of other elements. The analytical method to involve the crack in the structure corresponding to the second technique can be found in, e.g. [21]. The main motivation in this research is to investigate the results of the two techniques and to express the priority or the weakness of those. Furthermore, an investigation is performed to represent the advantages and drawbacks of the semi-analytical finite element method when different cracks are incorporated within the cylinder. A buckling analysis is carried out to compare the results of two techniques. Also, the effects of the geometry of the cylinder, the crack depth and the crack position on the critical buckling load are evaluated. Moreover, the vibration analysis is conducted to investigate the influence of the circumferential crack on the natural frequency.

## 2. Model

An isotropic thin cylindrical shell with a circumferential crack is considered in the analysis as shown in Fig. 1. The cylinder is modeled by the one-dimensional model in the framework of the semi-analytical finite element method.

Moreover, the crack is modeled by using a rotational spring corresponding to geometric and material characteristics of the crack.

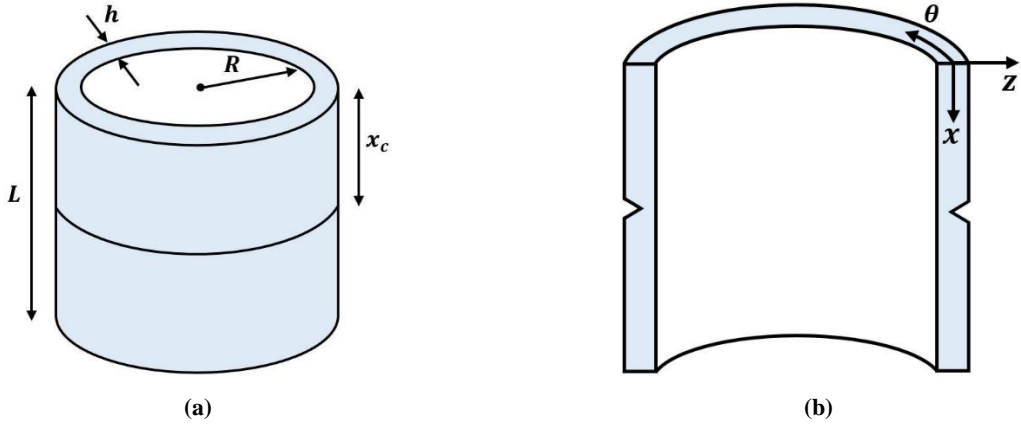

**Fig. 1** Cracked cylindrical shell: **a)** geometric parameters; **b)** longitudinal cross-section.

### 2.1. Modeling of the crack

Fig. 2 shows a cylindrical shell including the crack modeled by a rotational spring. Many research works on beams and cylinders have been carried out to present relations between the stiffness factor of the rotational spring and the crack characteristics, see e.g. [13] and [22]. The investigations show that different relations presented in the references to determine the stiffness factor of the rotational spring result in similar outputs.

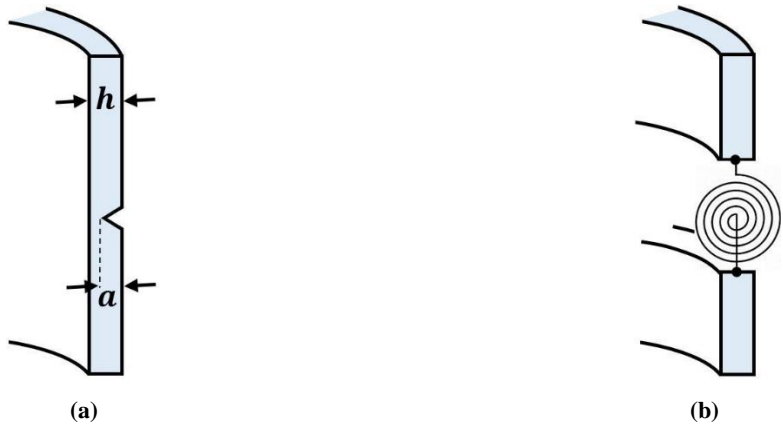

**Fig. 2** Crack in the cylinder: **a)** circumferential cracked cylinder; **b)** rotational spring model

Accordingly, equations presented by Yokoyama and Chen [13], which have already been used in the analysis of a cracked Euler-Bernoulli beam, are employed to model the crack as a circumferentially distributed rotational spring for the cylindrical shell by considering $b = 2\pi R$ and $\mu = \frac{a}{h}$.

$$K_I = \frac{6M_x}{bh^2}\sqrt{\pi a}F_M(\xi) \qquad for\ 0 \leq \mu \leq 0.6 \tag{1a}$$

$$K_I = \frac{3.99 M_x}{bh\sqrt{h}\sqrt{(1-\mu)^3}} \quad for\ 0.6 < \mu < 1.0 \tag{1b}$$

in which $M_x$ is the bending moment and

$$F_M(\mu) = \sqrt{\left(\frac{2}{\pi\mu}\right)\tan\frac{\pi\mu}{2}}\frac{0.923 + 0.199[1 - \sin\left(\frac{\pi\mu}{2}\right)]^4}{\cos\left(\frac{\pi\mu}{2}\right)} \tag{1c}$$

The range of $\mu = \frac{a}{h}$ for the research applications has been given in [13] and [21]. The stiffness factor of the spring is obtained as follows

$$\frac{1}{k_s} = \frac{2b(1-\nu^2)}{E}\int_0^a \left(\frac{K_I}{M_x}\right)^2 da \tag{2}$$

### 2.2. Modeling of cylinder

A one-dimensional model based on the semi-analytical finite element method is used in the analysis. The semi-analytical finite element method is formulated by one-dimensional elements in the axial direction, the first-order shear deformation theory in the radial direction and Fourier series in the circumferential direction [29-31].

The modal deformation can be formulated as follows

$$\begin{aligned} u(x,\theta) &= u_e\cos(n\theta) \\ v(x,\theta) &= v_e\sin(n\theta) \\ w(x,\theta) &= w_e\cos(n\theta) \\ \phi(x,\theta) &= \phi_e\cos(n\theta) \end{aligned} \tag{3}$$

The governing equations and the description of the terms in Eq. (3) can be found in Appendix A. Uncracked elements, so-called standard elements, are assumed to have two nodes and eight degrees of freedom as shown in Fig. 3.

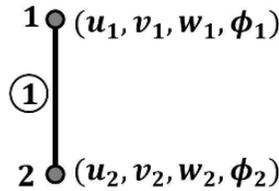

**Fig. 3** Degrees of freedom for an intact element

The displacement field, shape functions and isoparametric formulation are described in Appendix B. The kinematic equations are established by the Kirchhoff hypotheses in the theory of cylindrical shells [32] as

$$\boldsymbol{\varepsilon}_L = \begin{Bmatrix} \varepsilon_{x0} \\ \varepsilon_{\theta 0} \\ \gamma_{x\theta 0} \\ \kappa_x \\ \kappa_\theta \\ \kappa_{x\theta} \end{Bmatrix} = \begin{Bmatrix} \dfrac{\partial u}{\partial x} \\ \dfrac{1}{R}\left(\dfrac{\partial v}{\partial \theta} - w\right) \\ \dfrac{1}{R}\dfrac{\partial u}{\partial \theta} + \dfrac{\partial v}{\partial x} \\ \dfrac{\partial^2 w}{\partial x^2} \\ \dfrac{1}{R^2}\left(\dfrac{\partial v}{\partial \theta} + \dfrac{\partial^2 w}{\partial \theta^2}\right) \\ \dfrac{1}{R}\left(\dfrac{\partial v}{\partial x} + \dfrac{\partial^2 w}{\partial x \partial \theta}\right) \end{Bmatrix} \qquad (4)$$

The desired form of the strain in the finite element analysis is when derivatives of the shape functions and the nodal displacement are separated as

$$\begin{Bmatrix} \varepsilon_{x0} \\ \varepsilon_{\theta 0} \\ \gamma_{x\theta 0} \\ \kappa_x \\ \kappa_\theta \\ \kappa_{x\theta} \end{Bmatrix} = \begin{Bmatrix} \boldsymbol{B}_x \cos(n\theta)\vec{u} \\ \boldsymbol{B}_\theta \cos(n\theta)\vec{u} \\ \boldsymbol{B}_{x\theta} \sin(n\theta)\vec{u} \\ \boldsymbol{B}_{\kappa x} \cos(n\theta)\vec{u} \\ \boldsymbol{B}_{\kappa \theta} \cos(n\theta)\vec{u} \\ \boldsymbol{B}_{\kappa x\theta} \sin(n\theta)\vec{u} \end{Bmatrix} \qquad (5)$$

which components of $\boldsymbol{B}$ can be found in Appendix C.

The first-order shear deformation theory is utilized to relate the strains of the neutral axis and those of other points as

$$\boldsymbol{\varepsilon} = \begin{Bmatrix} \varepsilon_x \\ \varepsilon_\theta \\ \gamma_{x\theta} \end{Bmatrix} = \begin{Bmatrix} \varepsilon_{x0} - z\kappa_x \\ \varepsilon_{\theta 0} - z\kappa_\theta \\ \gamma_{x\theta 0} - 2z\kappa_{x\theta} \end{Bmatrix} \qquad (6)$$

The constitutive equation in a thin cylinder is given as

$$\boldsymbol{\sigma} = \begin{Bmatrix} \sigma_x \\ \sigma_\theta \\ \tau_{x\theta} \end{Bmatrix} = \begin{bmatrix} \dfrac{E}{1-\nu^2} & \dfrac{\nu E}{1-\nu^2} & 0 \\ \dfrac{\nu E}{1-\nu^2} & \dfrac{E}{1-\nu^2} & 0 \\ 0 & 0 & G \end{bmatrix} \begin{Bmatrix} \varepsilon_x \\ \varepsilon_\theta \\ \gamma_{x\theta} \end{Bmatrix} \qquad (7)$$

The strain energy of an element can be obtained as follows

$$U = \frac{1}{2}\int \boldsymbol{\sigma} \cdot \boldsymbol{\varepsilon}\, dV \qquad (8a)$$

whose stiffness matrix is extracted from the strain energy in the following form

$$U = \frac{1}{2}\boldsymbol{u}^T \boldsymbol{k}\boldsymbol{u} \qquad (8b)$$

in which the stiffness matrix is obtained as

$$\boldsymbol{k} = \int_{-1}^{1} \boldsymbol{B}^T \boldsymbol{D} \boldsymbol{B} \pi R \frac{l}{2} d\xi \qquad (9)$$

Moreover, the kinetic energy of an element can be computed as

$$KE = \frac{1}{2}\rho \int (\dot{u}^2 + \dot{v}^2 + \dot{w}^2) dV \tag{10a}$$

in which $\rho$ is the density. The kinetic energy can be written in the form of

$$KE = \frac{1}{2}\dot{u}^T m \dot{u} \tag{10b}$$

The mass matrix is derived using

$$m = \frac{\rho}{2} \int_{-1}^{1} N^T N \pi R \frac{l}{2} h d\xi \tag{11}$$

$D$, $B$ and $N$ are given in Appendix D.

## 2.3. Buckling analysis

An eigenvalue solution is applied in the linear buckling analysis in which two main parameters are the assembled stiffness matrix and geometric stiffness matrix.

$$|K + \lambda K_G| = 0 \tag{12}$$

The global geometric stiffness matrix, $K_G$, is obtained by considering the nonlinear strain and initial stress in the structure. The nonlinear strain is given as

$$\varepsilon_{NL} = \begin{Bmatrix} \frac{1}{2}\left(\frac{\partial w}{\partial x}\right)^2 \\ \frac{1}{2}\left(\frac{\partial w}{R \partial \theta}\right)^2 \\ \frac{1}{R}\frac{\partial w}{\partial x}\frac{\partial w}{\partial \theta} \end{Bmatrix} \tag{13}$$

The stress tensor corresponding to Eq. (7), which is used in the determination of the geometric stiffness matrix, is written as

$$\sigma = \begin{bmatrix} \sigma_x & \tau_{x\theta} \\ \tau_{x\theta} & \sigma_\theta \end{bmatrix} \tag{14}$$

Considering Eqs. (13) and (14), the geometric stiffness matrix of an element is given by

$$k_G = \int G^T \sigma G dV \tag{15}$$

which $G$ can be found in detail by [19], [29] and [33].

## 2.4. Vibration analysis

The vibration analysis is used to determine the natural frequency of the cylinder. An eigenvalue solution to specify the natural frequency is as follows

$$|K - M\omega^2| = 0 \tag{16}$$

in which $M$ is the global mass matrix. It is assumed that no change is made in components of the matrix due to the crack. In other words, the mass matrix is independent from the crack effect.

## 3. Description of techniques

Two techniques are implemented to quantify the crack effects in the analysis. The first one, which is originally applied for a cylindrical shell, is formulated based on the conversion matrix technique. On the other side, the second technique is introduced based on the definition of a stiffness matrix at the crack point through a set of springs equaled with the crack parameters. The crack parameters are involved in the global stiffness matrix when the stiffness matrix of the set of springs and stiffness matrices of the standard elements of the cylinder are assembled.

### 3.1. Conversion technique

This technique, which was already applied over cracked beams [16-18], is implemented by dividing a cracked element into three parts including two sub-elements and a rotational spring. The finite element method is used to obtain the enriched stiffness matrix from the strain energy of the three parts. Therefore, a cracked element as shown in Fig. 4 includes four nodes in which displacements of two middle nodes are obtained in terms of the two other nodes by considering continuity conditions.

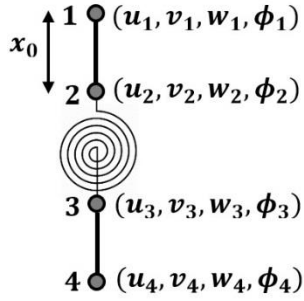

**Fig. 4** Sub-elements and degrees of freedom for a cracked element in conversion matrix technique

The boundary conditions should be satisfied at the crack point in which displacements and loads in two sides of the crack point are related to each other as

$$u_2 = u_3 \tag{17a}$$
$$v_2 = v_3 \tag{17b}$$
$$w_2 = w_3 \tag{17c}$$
$$N_{xT}(x_0) = N_{xB}(0) \tag{17d}$$
$$N_{x\theta T}(x_0) = N_{x\theta B}(0) \tag{17e}$$
$$M_{xT}(x_0) = M_{xB}(0) \tag{17f}$$
$$\phi_2 + \phi_s = \phi_3 \tag{17g}$$
$$Q_{xT}(x_0) = Q_{xB}(0) \tag{17h}$$

where forces and moments in the top and bottom sub-elements are denoted with *T* and *B* subscripts, respectively, in which [32]

$$N_x(u(x,\theta), v(x,\theta), w(x,\theta)) = \frac{Eh}{1-v^2}\left(\frac{\partial u}{\partial x} + \frac{v}{R}\left(\frac{\partial v}{\partial \theta} - w\right)\right) \tag{18a}$$

$$N_\theta(u(x,\theta), v(x,\theta), w(x,\theta)) = \frac{Eh}{1-v^2}\left(\frac{1}{R}\left(\frac{\partial v}{\partial \theta} - w\right) + v\frac{\partial u}{\partial x}\right)$$

$$N_{x\theta}(u(x,\theta), v(x,\theta)) = \frac{Eh}{2(1+v)}\left(\frac{\partial v}{\partial x} + \frac{1}{R}\frac{\partial u}{\partial \theta}\right)$$

$$M_x(v(x,\theta), w(x,\theta)) = -D\left(\frac{\partial^2 w}{\partial x^2} + \frac{v}{R^2}\left(\frac{\partial v}{\partial \theta} + \frac{\partial^2 w}{\partial \theta^2}\right)\right)$$

$$M_\theta(v(x,\theta), w(x,\theta)) = -D\left(\frac{1}{R^2}\left(\frac{\partial v}{\partial \theta} + \frac{\partial^2 w}{\partial \theta^2}\right) + v\frac{\partial^2 w}{\partial x^2}\right)$$

$$M_{x\theta}(v(x,\theta), w(x,\theta)) = -D\frac{(1-v)}{R}\left(\frac{\partial v}{\partial x} + \frac{\partial^2 w}{\partial x \partial \theta}\right)$$

$$Q_x(v(x,\theta), w(x,\theta)) = \frac{\partial M_x}{\partial x}$$

and

$$D = \frac{Eh^3}{12(1-v^2)} \tag{18b}$$

The relation of $Q_x$ is written by neglecting the effect of the torsional moment, $M_{x\theta}$, in the shear force. Two other continuity equations (i.e., $N_{\theta T}(x_0) = N_{\theta B}(0)$ and $M_{\theta T}(x_0) = M_{\theta B}(0)$) are dependent equations that yield similar relations to Eq. (17). Applying these boundary conditions gives two conversion matrices that are used to derive the stiffness matrix of the cracked element. Eight independent continuity conditions, Eq. (17), are applied to determine the displacements of the middle nodes $(u_2, v_2, w_2, \phi_2, u_3, v_3, w_3, \phi_3)$ with respect to displacements of the main nodes $(u_1, v_1, w_1, \phi_1, u_4, v_4, w_4, \phi_4)$. It can be represented as follows

$$\vec{u}_T = \boldsymbol{C}_T \vec{u} \tag{19a}$$

$$\vec{u}_B = \boldsymbol{C}_B \vec{u} \tag{19b}$$

in which $\vec{u}_T$ and $\vec{u}_B$ are denoted as the displacement vector of top and bottom sides sub-elements, respectively, and $\vec{u}$ is the displacement vector of the cracked element. These vectors are defined as $\vec{u}_T = [u_1, v_1, w_1, \phi_1, u_2, v_2, w_2, \phi_2]^T$, $\vec{u}_B = [u_3, v_3, w_3, \phi_3, u_4, v_4, w_4, \phi_4]^T$ and $\vec{u} = [u_1, v_1, w_1, \phi_1, u_4, v_4, w_4, \phi_4]^T$. $\boldsymbol{C}_T$ and $\boldsymbol{C}_B$ are called conversion matrices related to top and bottom sides sub-elements, respectively, described in Appendix E.

The conversion matrix technique is an energy-based technique in which the stiffness matrix of a cracked element is obtained by strain energies of two sub-elements and the rotational spring. This stiffness matrix is enriched through crack characteristics equaled in the spring. The sum of the strain energies in the cracked element is

$$U = U_T + U_B + U_{sp} \tag{20a}$$

Considering Eq. (8b) yields

$$\frac{1}{2}\vec{u}^T \boldsymbol{k}_{cr} \vec{u} = \frac{1}{2}\vec{u}_T^T \boldsymbol{k}_T \vec{u}_T + \frac{1}{2}\vec{u}_B^T \boldsymbol{k}_B \vec{u}_B + \frac{1}{2}k_s(\phi_3 - \phi_2)^2 \tag{20b}$$

Rotations of $\phi_2$ and $\phi_3$ can be written in terms of displacements of main nodes as

$$\phi_2 = \boldsymbol{C}_{T\phi_2}\vec{u} \tag{21a}$$

$$\phi_3 = \boldsymbol{C}_{B\phi_3}\vec{u} \tag{21b}$$

in which $C_{T\phi_2}$ and $C_{B\phi_3}$ are eighth row of the top conversion matrix and fourth row of the bottom one, respectively. Therefore, the stiffness matrix of a cracked element is determined by the substitution of Eqs. (19) and (21) into (20b)

$$\boldsymbol{k}_{cr} = \boldsymbol{C}_T^T \boldsymbol{k}_T \boldsymbol{C}_T + \boldsymbol{C}_B^T \boldsymbol{k}_B \boldsymbol{C}_B + k_s \big(\boldsymbol{C}_{B\phi_3} - \boldsymbol{C}_{T\phi_2}\big)^T \big(\boldsymbol{C}_{B\phi_3} - \boldsymbol{C}_{T\phi_2}\big) \tag{22}$$

Eq. (9) is used to obtain the top and bottom stiffness matrices (i.e., $\boldsymbol{k}_T$ and $\boldsymbol{k}_B$). An analogous way results in the determination of the geometric stiffness matrix

$$\boldsymbol{k}_{Gcr} = \boldsymbol{C}_T^T \boldsymbol{k}_{GT} \boldsymbol{C}_T + \boldsymbol{C}_B^T \boldsymbol{k}_{GB} \boldsymbol{C}_B \tag{23}$$

where the top and bottom geometric stiffness matrices, $\boldsymbol{k}_{GT}$ and $\boldsymbol{k}_{GB}$, are determined based on Eq. (15) in which the conversion matrices are applied to calculate stress tensors of top and bottom sub-elements.

### 3.2. Spring set technique

In this technique, a separate element as a set of springs is considered to involve the crack parameters into the global stiffness matrix. In other words, the global stiffness matrix of the cracked cylinder is obtained without considering any sub-elements, unlike the conversion matrix technique.

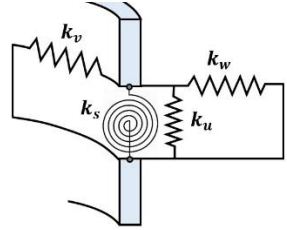

**Fig. 5** Spring set instead of crack

Fig. 5 shows four springs to model the crack effect on the stiffness structure. Three springs represented in the axial, radial and circumferential directions are used to satisfy the continuity conditions in the three mentioned directions. The stiffness factors of the three springs are taken into account considerable amounts to apply the geometric boundary conditions at the crack point. The rotational spring explains the crack characteristics as obtained from Eq. (2). The stiffness matrix related to these four springs is given by

$$\boldsymbol{K}_{\text{crack}} = \begin{bmatrix} k_u & 0 & 0 & 0 & -k_u & 0 & 0 & 0 \\ 0 & k_v & 0 & 0 & 0 & -k_v & 0 & 0 \\ 0 & 0 & k_w & 0 & 0 & 0 & -k_w & 0 \\ 0 & 0 & 0 & k_s & 0 & 0 & 0 & -k_s \\ -k_u & 0 & 0 & 0 & k_u & 0 & 0 & 0 \\ 0 & -k_v & 0 & 0 & 0 & k_v & 0 & 0 \\ 0 & 0 & -k_w & 0 & 0 & 0 & k_w & 0 \\ 0 & 0 & 0 & -k_s & 0 & 0 & 0 & k_s \end{bmatrix} \tag{24}$$

in which the stiffness matrix obtained at the crack point as a complete element is assembled with other elements of the cylinder. Fig. 6 shows three elements for a section of the cylindrical shell, which the first and third elements are considered as one-dimensional standard elements related to the semi-analytical finite element method, while the second element at the crack point has been added to the structure to represent the softness due to the crack. The stiffness matrix of the second element is introduced as Eq. (24).

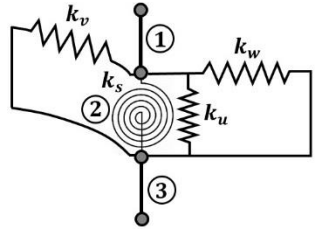

**Fig. 6** Schematic description of two standard elements and one cracked element

## 4. Results and discussion

The modal deformation of some case studies related to buckling and vibration analyses is examined to determine the results of the critical buckling load and the natural frequency using different methods including two presented techniques, ABAQUS modeling, and experimental test.

### 4.1. The buckling of the cracked cylinder

The validation of techniques represented for the cracked cylindrical shell is evaluated through a case study mentioned in [21] with geometric and material characteristics as $(h = 0.2m, m = 1)$ and $(\nu = 0.3, E = 200\text{GPa})$, respectively, and simply supported end conditions in which $m = \frac{12(1-\nu^2)}{R^2h^2}$ and the length of the cylinder is selected large enough with the circumferential crack in the middle. Table 1 compares the results of the two techniques with [21] in which the two techniques give close outputs to each other and the reference. The critical buckling loads mentioned in Table 1 are related to the first circumferential mode. As it is seen from Table 1, increasing the crack depth results in decreasing the critical buckling load. The most drastic decrease is related to $\frac{a}{h} = 0.9$ in which the critical load of the cracked cylinder is approximately half of the critical load of the intact cylinder. A nonlinear behavior is observed between the crack depth and the critical buckling load, as the buckling load capacity in $\frac{a}{h} \leq 0.5$ decreases nearly 10%, while it reduces almost 50% for $0.5 \leq \frac{a}{h} \leq 0.9$.

**Table 1**

Validation of the two techniques: I) Conversion technique; II) Spring set technique; III) [21]

| $\frac{a}{h}$ | 0 | 0.1 | 0.2 | 0.3 | 0.4 | 0.5 | 0.6 | 0.7 | 0.8 | 0.9 |
|---|---|---|---|---|---|---|---|---|---|---|
| I | 2.00 | 2.00 | 1.99 | 1.97 | 1.94 | 1.87 | 1.74 | 1.50 | 1.27 | 1.08 |
| II | 2.00 | 2.00 | 1.99 | 1.97 | 1.94 | 1.87 | 1.73 | 1.49 | 1.25 | 1.07 |
| III | 2.00 | 2.00 | 1.99 | 1.98 | 1.95 | 1.88 | 1.73 | 1.50 | 1.26 | 1.07 |

The convergence of results of the two techniques implemented into the framework of the semi-analytical finite element method is investigated by Figs. 7a, 7b and 7c for the conversion and set spring techniques. Fig. 7 shows that the discretization of the cylinder via 21 one-dimensional elements yields acceptable results. In other words, the comparison of the first two curves in 21 elements to 41 elements confirms that the number of 21 elements is an appropriate selection to analyze, However a close agreement is seen between different elements and also two techniques in Fig. 7c. Material and geometric properties are considered like what were mentioned in Table 1 with $\frac{a}{h} = 0.5$.

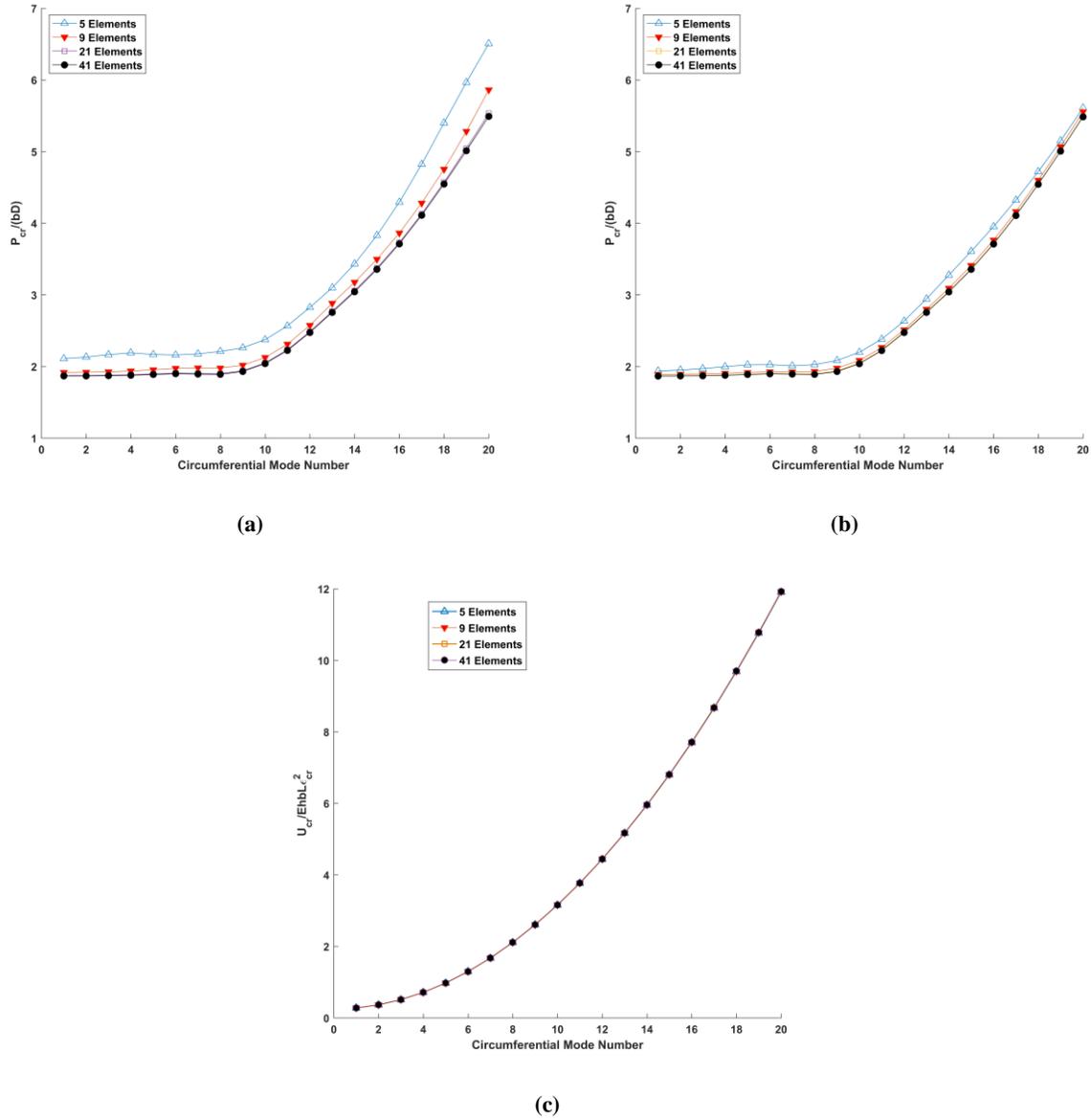

**Fig. 7** Investigation of convergence: **a)** Critical buckling load in conversion technique; **b)** Critical buckling load in spring set technique; **c)** Critical strain energy in both techniques

The deformed shape of the cylinder under the critical buckling load has been displayed in Fig. 8 for different circumferential mode numbers by inserting Matlab results into Tecplot software. The deformed shape is considered to be ten times of the real value for the clear visibility of mode shapes. The circumferential crack is assumed at the middle of the cylinder with $\frac{a}{h} = 0.5$, whose effect is evaluated by decreasing the stiffness in the cracked element without any change in the appearance of geometry.

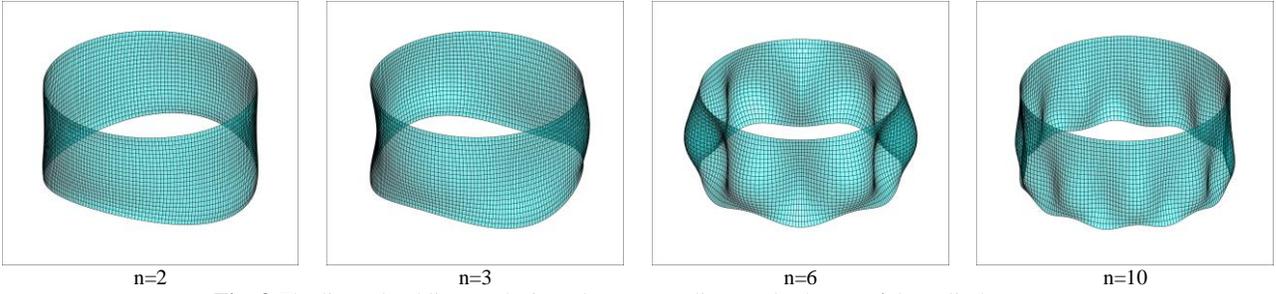

| n=2 | n=3 | n=6 | n=10 |

**Fig. 8** The linear buckling analysis and corresponding mode shapes of the cylinder

An experimental study has been also carried out to investigate the influence of the circumferential crack in the critical buckling load. Fig. 9 shows a cylindrical shell under the axial compression in a uniaxial compressive test. The length, radius, and thickness of the cylinder are considered 100, 115 and 1 mm, respectively, with material properties mentioned in Table 1.

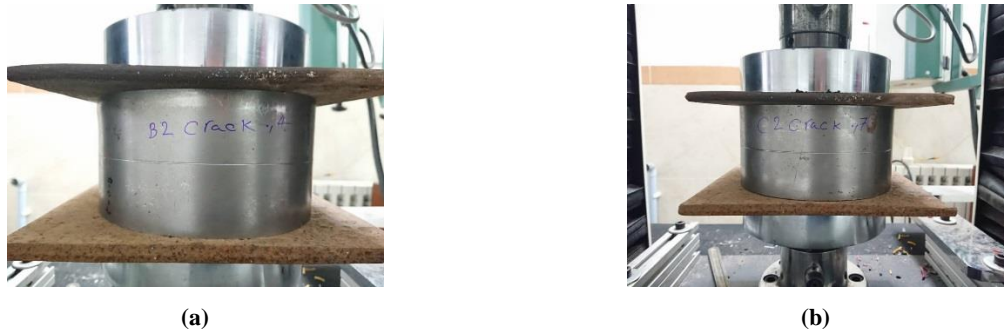

(a)         (b)

**Fig. 9** Cracked cylindrical shell in the uniaxial compressive test: **a)** $\frac{a}{h} = 0.4$ ; **b)** $\frac{a}{h} = 0.7$

Table 2 shows a comparison between the results of the experimental test, semi-analytical finite element method and ABAQUS modeling in which the critical buckling load decreased due to the initial circumferential crack on the cylinder has been determined.

**Table 2**

Reduction of the critical buckling load. I) Experimental test; II) Conversion technique; III) ABAQUS

| $\frac{a}{h}$ | $\frac{P_{\text{cracked}}}{P_{\text{intact}}}$ | | |
|---|---|---|---|
| | I | II | III |
| 0.4 | 0.84 | 0.98 | 0.98 |
| 0.7 | 0.67 | 0.77 | 0.70 |

Table 2 demonstrates the result of ABAQUS is between the conversion technique and experimental results for a/h=0.7 and it is equal to the value of the conversion technique for a/h=0.4. Beside the comparison of results of three methods in Table. 2, the main aim to apply ABAQUS in the analysis is to represent advantages of the presented finite element method in detailed and quantitatively. Some significant disadvantages of ABAQUS observed in this modeling are listed as follows:

ABAQUS sorts the results of the buckling analysis just in terms of eigenvalues. Moreover, it gives a mixture of buckling modes without the separation of the circumferential and axial modes.

Besides the global partitioning of the cylinder in ABAQUS, the crack zone should be partitioned for each the crack depth, as the change of the depth leads to re-partitioning of the crack zone.

The convergence of ABAQUS results requires much more time consuming than the convergence of the presented techniques results. Table 3 represents a quantitative comparison between the convergence of the two theoretical methods to determine the critical buckling load as

**Table 3**

Comparison of the convergence between the conversion technique and ABAQUS

| Method | Mesh size(mm) | Time (s) | $\frac{P}{P_{conv}}$ | Mesh size(mm) | Time (s) | $\frac{P}{P_{conv}}$ | Mesh size(mm) | Time (s) | $\frac{P}{P_{conv}}$ | Mesh size(mm) | Time (s) | $\frac{P}{P_{conv}}$ |
|---|---|---|---|---|---|---|---|---|---|---|---|---|
| Conversion technique (1D) | 20 | 2.3 | 1.09 | 9 | 2.9 | 1.01 | 5 | 4.4 | 1.001 | 2.5 | 12.6 | 1 |
| ABAQUS (3D) | $10 \times 10$ | 25 | 2.07 | $8 \times 8$ | 210 | 1.47 | $2 \times 2$ | 2610 | 1.003 | $1 \times 1$ | 7318 | 1 |

The data of Table 3 produced via a usual personal computer are related to the intact cylinder with geometric and material properties mentioned in Table 2. It is obvious that results obtained from presented technique give faster convergence in comparison with ABAQUS. The main reason of this good convergence is to combine one-dimensional model with the analytical method called the semi-analytical finite element method.

Table 4 represents the influence of the crack position on the critical buckling load with the characteristics similar to Table 1 and ($\frac{a}{h} = 0.5, L = 5\pi$). Results show that if the crack sits at $\frac{x_c}{L} = 0.2$ or 0.8, the softening of the structure due to the crack can be ignored. Also, the crack around the edge conditions (i.e. $\frac{x_c}{L} = 0.1$) leads to a maximum decrease in the critical buckling load.

**Table 4**

Effect of the crack position on the critical buckling load

| $\frac{x_c}{L}$ | 0.1 | 0.2 | 0.3 | 0.4 | 0.5 |
|---|---|---|---|---|---|
| $\frac{P_{cr}}{bD}$ | 1.80 | 2.00 | 1.86 | 1.95 | 1.87 |

The effect of the crack depth on the critical load in different circumferential mode numbers is investigated in Fig. 10 with the characteristics considered in Table 1. Fig. 10 displays that there is a nonlinear behavior between the crack depth and the critical buckling load. Also, the critical buckling load is approximately constant for the circumferential mode number less than 9, while an increasable nonlinear behavior is seen in the critical buckling curves when the circumferential mode number is more than 8.

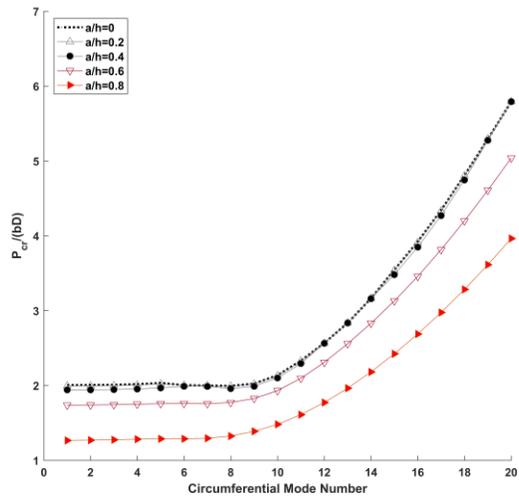

**Fig. 10** Effect of the crack depth on the critical buckling load

Fig. 11 demonstrates the influence of the length of the cracked cylindrical shell on the critical buckling load. The curves show when the length of the cylinder increases, results are converged to a permanent state. In other words, the critical buckling load can be considered independent of the cylinder length after a certain one, e.g. $L = 3\pi$.

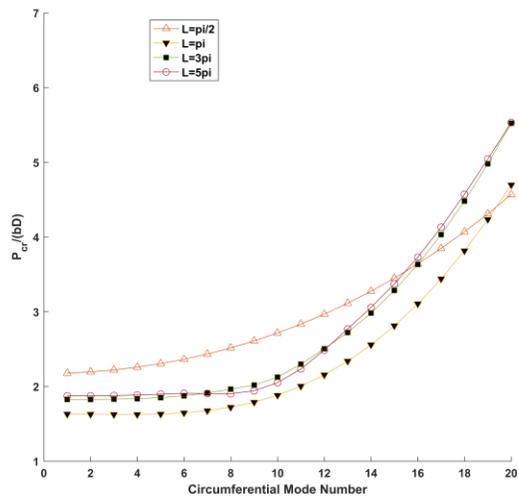

**Fig. 11** Effect of the length on the critical buckling load

Results and mentioned relations show that one of the main advantages of the two techniques in comparison with previous research works is to involve the circumferential mode in the analysis. the semi-analytical finite element in combination with the two techniques makes a powerful and efficient procedure in modal analyses. This procedure can be finely implemented in the nonlinear analysis like the nonlinear vibration or post-buckling problems. In other words, the procedure of the nonlinear analysis for cracked cylindrical shells can be implemented quickly and at a low cost in which the cost of the represented techniques is less than the cost of the general-purposes programs. On the other side, the represented conversion technique can be efficiently employed in the modal nonlinear analysis or circumferentially asymmetric cracks, where analytical methods may be incapable of the solution of such problems.

The essential difference between the two mentioned techniques (i.e. conversion and spring set) is related to the cracked element. The conversion technique introduces an enriched stiffness matrix in the cracked element without adding the number of degrees of freedom, while a set of springs instead of the crack used in the spring set technique leads to the increase of the number of degrees of freedom of the structure. Results obtained from the two techniques show a good agreement, however obvious differences in relations are observed. In other words, the conversion technique has been implemented in a more comprehensive framework than the spring set technique, while its established procedure is more complicated.

### 4.2. The vibration of the cracked cylinder

The validity of the derived equations of the natural frequency based on the semi-analytical finite element method is evaluated for an intact cylindrical shell. Fig. 12 demonstrates the effect of the circumferential mode number on the frequency parameter in which results show quite close to [34]. The cylinder is considered the simply supported-simply supported ($v = w = M_x = N_x = 0$) with *R/h*=500.

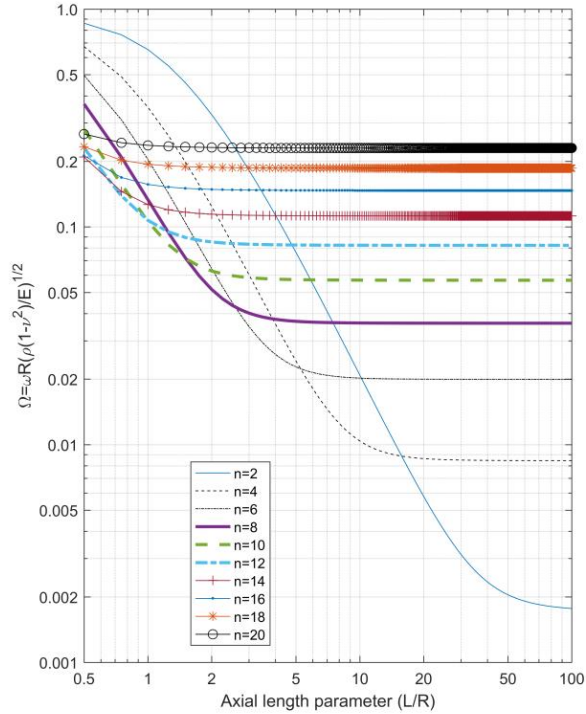

**Fig. 12** Effect of the mode number on the frequency parameter of intact cylinder

The effect of the crack is investigated by involving the enriched stiffness matrix into the eigenvalue equation of the vibration. The density ($\rho$) is assumed to be 7850 kg/m$^3$ with the boundary conditions similar to Fig. 12 and the geometric and material characteristics mentioned in Table 1. The results of Table 5 show that the effect of the circumferential crack in the natural frequency of the cylinder is negligible. An interesting comment that had been already mentioned in [35].

**Table 5**

Effect of the crack depth and mode number on the frequency parameter ($\Omega = \omega R \sqrt{\frac{\rho(1-\nu^2)}{E}}$)

|  |  | \multicolumn{4}{c}{$a/h$} |
| --- | --- | --- | --- | --- | --- |
|  |  | 0.1 | 0.3 | 0.5 | 0.8 |
| Mode Number | 1 | 0.8557 | 0.8557 | 0.8557 | 0.8557 |
|  | 3 | 0.5172 | 0.5172 | 0.5172 | 0.5169 |
|  | 7 | 0.2899 | 0.2896 | 0.2890 | 0.2867 |
|  | 11 | 0.4681 | 0.4672 | 0.4662 | 0.4620 |

## 5. Conclusion

In this paper, two techniques, called *conversion* and *spring set*, into the framework of the semi-analytical finite element method have been initially introduced to determine the modal deformation of the cylindrical shells including a circumferential crack. An experimental test has been also carried out to investigate the reduction of the buckling load of the cracked cylindrical shell. The effect of the circumferential mode number on the critical buckling load of the cracked cylinder has been originally investigated. The validity of relations of the two techniques has been evaluated by the results of references, the experimental test, and ABAQUS modeling. One of the essential advantages of these two techniques especially the conversion technique is the feasibility of the development for the nonlinear or large deformation or asymmetric problems. These techniques are specifically effective and problem-solving in the nonlinear modal analysis where the analytical solutions or the empirical study or commercial personal computer programs may be incapable or too costly to analyze circumferentially cracked cylindrical shells. Results show that the maximum crack depth can decrease half of the load-bearing capacity, and a nonlinear behavior is observed between the crack depth and the critical buckling load. Moreover, the influence of the circumferential crack on the natural frequency can be ignored.

**Appendix A**

The relatively simple Donnell type shell theory can be used to analyze the shell stability in which the differential equations of the equilibrium is approximately written in the following form [36] and [37]

$$\frac{\partial N_x}{\partial x} + \frac{\partial N_{x\theta}}{R\partial \theta} = 0 \tag{A.1}$$

$$\frac{\partial N_{x\theta}}{\partial x} + \frac{\partial N_\theta}{R\partial \theta} = 0 \tag{A.2}$$

$$\frac{\partial^2 M_x}{\partial x^2} + \frac{2\partial^2 M_{x\theta}}{R\partial x\partial \theta} + \frac{\partial^2 M_\theta}{R^2\partial^2 \theta} - \frac{N_\theta}{R} + N_x \frac{\partial^2 w}{\partial x^2} + 2N_{x\theta} \frac{\partial^2 w}{R\partial x\partial \theta} + N_\theta \frac{\partial^2 w}{R^2\partial^2 \theta} = 0 \tag{A.3}$$

The Airy stress function $f$ is utilized in the analysis as

$$N_x = \frac{\partial^2 f}{R^2 \partial \theta^2}, N_\theta = \frac{\partial^2 f}{\partial x^2}, N_{x\theta} = -\frac{\partial^2 f}{R \partial x \partial \theta} \tag{A.4}$$

A small perturbation is used to trace the post-buckling path via

$$w \coloneqq w_0 + w \tag{A.5}$$

$$f \coloneqq f_0 + f \tag{A.6}$$

in which subscript 0 denotes parameters in initial state of the cylinder. Therefore, the stability equation is obtained by the substitution of Eqs. (18a), (A.4) and (A.5) into Eq. (A.3). The variables separable form is used in the buckling solution of the cylinder based on trigonometric functions of Eq. (3).

**Appendix B**

The displacement field in different directions for an element is obtained in terms of shape functions and the nodal displacements as

$$\begin{aligned} u_e &= N_1 u_1 + N_2 u_2 \\ v_e &= N_1 v_1 + N_2 v_2 \\ w_e &= H_1 w_1 + H_2 \phi_1 + H_3 w_2 + H_4 \phi_2 \\ \phi_e &= \frac{\partial w_e}{\partial x} = \frac{\partial H_1}{\partial x} w_1 + \frac{\partial H_2}{\partial x} \phi_1 + \frac{\partial H_3}{\partial x} w_2 + \frac{\partial H_4}{\partial x} \phi_2 \end{aligned} \tag{B.1}$$

in which Lagrange and Hermite shape functions are utilized to interpolate $(u, v)$ and $(w, \phi)$, respectively

$$N_1(x) = 1 - \frac{x}{l}, N_2(x) = \frac{x}{l} \tag{B.2}$$

$$H_1(x) = \frac{1}{l^3}(2x^3 - 3x^2 l + l^3), H_2(x) = \frac{1}{l^3}(x^3 L - 2x^2 l^2 + x l^3) \tag{B.3}$$

$$H_3(x) = \frac{1}{l^3}(-2x^3 + 3x^2 l), H_4(x) = \frac{1}{l^3}(x^3 l - x^2 l^2)$$

The Gauss integration is applied by using the local coordinate $\xi$ $[-1 \leq \xi \leq 1]$, which is defined as

$$x(\xi) = \frac{l}{2}(1 + \xi) \quad \begin{bmatrix} \xi = -1 \\ x = 0 \end{bmatrix} \begin{bmatrix} \xi = 1 \\ x = l \end{bmatrix} \tag{B.4}$$

Therefore, the shape functions can be rewritten with respect to the local coordinate whereby the integral equations of the finite element method should be transformed as

$$\int_0^l g(x) dx = \int_{-1}^1 g(\xi) \frac{dx}{d\xi} d\xi = \int_{-1}^1 g(\xi) \frac{l}{2} d\xi = \sum_{i=1}^3 g(\xi_i) \frac{l}{2} W_i \tag{B.5}$$

The one-dimensional Gauss integration with three Gauss points [38] is used in the analysis in which the weighting factors and the coordinate of the Gauss points are $W_i$ and $\xi_i$, respectively.

**Appendix C**

$$\boldsymbol{B}_x = [N_1 \ 0 \ 0 \ 0 \ N_2 \ 0 \ 0 \ 0] \tag{C.1}$$

$$\boldsymbol{B}_\theta = \left[0 \ \frac{n}{R} N_1 \ -\frac{1}{R} H_1 \ -\frac{1}{R} H_2 \ 0 \ \frac{n}{R} N_2 \ -\frac{1}{R} H_3 \ -\frac{1}{R} H_4\right] \tag{C.2}$$

$$\boldsymbol{B}_{x\theta} = \left[-\frac{n}{R}N_1 \quad N_{1,x} \quad 0 \quad 0 \quad -\frac{n}{R}N_2 \quad N_{2,x} \quad 0 \quad 0\right] \tag{C.3}$$

$$\boldsymbol{B}_{\kappa x} = [0 \quad 0 \quad H_{1,xx} \quad H_{2,xx} \quad 0 \quad 0 \quad H_{3,xx} \quad H_{4,xx}] \tag{C.4}$$

$$\boldsymbol{B}_{\kappa\theta} = \left[0 \quad \frac{n}{R^2}N_1 \quad -\frac{n^2}{R^2}H_1 \quad -\frac{n^2}{R^2}H_2 \quad 0 \quad \frac{n}{R^2}N_2 \quad -\frac{n^2}{R^2}H_3 \quad -\frac{n^2}{R^2}H_4\right] \tag{C.5}$$

$$\boldsymbol{B}_{\kappa x\theta} = \left[0 \quad \frac{1}{R}N_{1,x} \quad -\frac{n}{R}H_{1,x} \quad -\frac{n}{R}H_{2,x} \quad 0 \quad \frac{1}{R}N_{2,x} \quad -\frac{n}{R}H_{3,x} \quad -\frac{n}{R}H_{4,x}\right] \tag{C.6}$$

**Appendix D**

$$\boldsymbol{B} = [\boldsymbol{B}_x^T \quad \boldsymbol{B}_\theta^T \quad \boldsymbol{B}_{x\theta}^T \quad \boldsymbol{B}_{\kappa x}^T \quad \boldsymbol{B}_{\kappa\theta}^T \quad \boldsymbol{B}_{\kappa x\theta}^T] \tag{D.1}$$

$$\boldsymbol{D} = \frac{12D}{h^2}\begin{bmatrix} 1 & \nu & 0 & 0 & 0 & 0 \\ \nu & 1 & 0 & 0 & 0 & 0 \\ 0 & 0 & \frac{h^2}{12} & \frac{\nu h^2}{12} & 0 & 0 \\ 0 & 0 & \frac{\nu h^2}{12} & \frac{h^2}{12} & 0 & 0 \\ 0 & 0 & 0 & 0 & \frac{1-\nu}{2} & 0 \\ 0 & 0 & 0 & 0 & 0 & \frac{(1-\nu)h^2}{6} \end{bmatrix} \tag{D.2}$$

$$\boldsymbol{N} = \begin{bmatrix} N_1 & 0 & 0 & 0 & N_2 & 0 & 0 & 0 \\ 0 & N_1 & 0 & 0 & 0 & N_2 & 0 & 0 \\ 0 & 0 & H_1 & H_2 & 0 & 0 & H_3 & H_4 \end{bmatrix} \tag{D.3}$$

**Appendix E**

$$\boldsymbol{C}_T = \begin{bmatrix} 1 & 0 & 0 & 0 & 0 & 0 & 0 & 0 \\ 0 & 1 & 0 & 0 & 0 & 0 & 0 & 0 \\ 0 & 0 & 1 & 0 & 0 & 0 & 0 & 0 \\ 0 & 0 & 0 & 1 & 0 & 0 & 0 & 0 \\ 1-\frac{x_0}{l} & 0 & 0 & 0 & \frac{x_0}{l} & 0 & 0 & 0 \\ 0 & 1-\frac{x_0}{l} & 0 & 0 & 0 & \frac{x_0}{l} & 0 & 0 \\ 0 & w_{v1} & w_{w1} & w_{f1} & 0 & w_{v4} & w_{w4} & w_{f4} \\ 0 & f2_{v1} & f2_{w1} & f2_{f1} & 0 & f2_{v4} & f2_{w4} & f2_{f4} \end{bmatrix} \tag{E.1}$$

$$C_B = \begin{bmatrix} 1-\dfrac{x_0}{l} & 0 & 0 & 0 & \dfrac{x_0}{l} & 0 & 0 & 0 \\ 0 & 1-\dfrac{x_0}{l} & 0 & 0 & 0 & \dfrac{x_0}{l} & 0 & 0 \\ 0 & w_{v1} & w_{w1} & w_{f1} & 0 & w_{v4} & w_{w4} & w_{f4} \\ 0 & f3_{v1} & f3_{w1} & f3_{f1} & 0 & f3_{v4} & f3_{w4} & f3_{f4} \\ 0 & 0 & 0 & 0 & 1 & 0 & 0 & 0 \\ 0 & 0 & 0 & 0 & 0 & 1 & 0 & 0 \\ 0 & 0 & 0 & 0 & 0 & 0 & 1 & 0 \\ 0 & 0 & 0 & 0 & 0 & 0 & 0 & 1 \end{bmatrix} \tag{E.2}$$

in which

$$w_{v1} = \frac{1}{\chi}(Dl^4n^3v^2x_0^3 - 4Dl^3n^3v^2x_0^4 + 6Dl^2n^3v^2x_0^5 - 4Dln^3v^2x_0^6 + Dn^3v^2x_0^7 - 6Dl^3R^2nvx_0^2 \tag{E.3}$$
$$+ 18Dl^2R^2nvx_0^3 - 18DlR^2nvx_0^4 + 6DR^2nvx_0^5)$$

$$w_{v4} = \frac{1}{\chi}(Dl^3n^3v^2x_0^4 - 3Dl^2n^3v^2x_0^5 + 3Dln^3v^2x_0^6 - Dn^3v^2x_0^7 - 6Dl^2R^2nvx_0^3 + 12DlR^2nvx_0^4 \tag{E.4}$$
$$- 6DR^2nvx_0^5)$$

$$w_{w1} = \frac{1}{\chi}(6Dl^3R^2n^2vx_0^2 - 18Dl^2R^2n^2vx_0^3 + 18DlR^2n^2vx_0^4 - 6DR^2n^2vx_0^5 + 3l^4R^4k_s - 9l^2R^4k_sx_0^2 \tag{E.5}$$
$$+ 6lR^4k_sx_0^3 + 12Dl^3R^4 - 36Dl^2R^4x_0 + 36DlR^4x_0^2 - 12DR^4x_0^3)$$

$$w_{w4} = \frac{1}{\chi}(6Dl^2R^2n^2vx_0^3 - 12DlR^2n^2vx_0^4 + 6DR^2n^2vx_0^5 + 9l^2R^4k_sx_0^2 - 6lR^4k_sx_0^3 + 12DR^4x_0^3) \tag{E.6}$$

$$w_{f1} = \frac{1}{\chi}(2Dl^3R^2n^2vx_0^3 - 6Dl^2R^2n^2vx_0^4 + 6DlR^2n^2vx_0^5 - 2DR^2n^2vx_0^6 + 3l^4R^4k_sx_0 - 6l^3R^4k_sx_0^2 \tag{E.7}$$
$$+ 3l^2R^4k_sx_0^3 + 12Dl^3R^4x_0 - 36Dl^2R^4x_0^2 + 36DlR^4x_0^3 - 12DR^4x_0^4)$$

$$w_{f4} = \frac{1}{\chi}(-2Dl^3R^2n^2vx_0^3 + 6Dl^2R^2n^2vx_0^4 - 6DlR^2n^2vx_0^5 + 2DR^2n^2vx_0^6 - 3l^3R^4k_sx_0^2 \tag{E.8}$$
$$+ 3l^2R^4k_sx_0^3 - 12DlR^4x_0^3 + 12DR^4x_0^4)$$

$$f2_{v1} = \frac{-1}{2\chi}(-3Dl^4n^3v^2x_0^2 + 15Dl^3n^3v^2x_0^3 - 27Dl^2n^3v^2x_0^4 + 21Dln^3v^2x_0^5 - 6Dn^3v^2x_0^6 \tag{E.9}$$
$$+ 24Dl^3R^2nvx_0 - 78Dl^2R^2nvx_0^2 + 90DlR^2nvx_0^3 - 36DR^2nvx_0^4)$$

$$f2_{v4} = \frac{-1}{2\chi}(-3Dl^3n^3v^2x_0^3 + 12Dl^2n^3v^2x_0^4 - 15Dln^3v^2x_0^5 + 6Dn^3v^2x_0^6 + 24Dl^2R^2nvx_0^2 \tag{E.10}$$
$$- 54DlR^2nvx_0^3 + 36DR^2nvx_0^4)$$

$$f2_{w1} = \frac{-1}{2\chi}(3Dl^4n^4v^2x_0^2 - 12Dl^3n^4v^2x_0^3 + 18Dl^2n^4v^2x_0^4 - 12Dln^4v^2x_0^5 + 3Dn^4v^2x_0^6 \tag{E.11}$$
$$- 24Dl^3R^2n^2vx_0 + 72Dl^2R^2n^2vx_0^2 - 72DlR^2n^2vx_0^3 + 24DR^2n^2vx_0^4$$
$$+ 36l^2R^4k_sx_0 - 36lR^4k_sx_0^2 + 36DR^4x_0^2)$$

$$f2_{w4} = \frac{-1}{2\chi}(-3Dl^2n^4v^2x_0^4 + 6Dln^4v^2x_0^5 - 3Dn^4v^2x_0^6 - 18Dl^2R^2n^2vx_0^2 + 36DlR^2n^2vx_0^3 \tag{E.12}$$
$$- 24DR^2n^2vx_0^4 - 36l^2R^4k_sx_0 + 36lR^4k_sx_0^2 - 36DR^4x_0^2)$$

$$f2_{f1} = \frac{-1}{2\chi}(Dl^4n^4v^2x_0^3 - 4Dl^3n^4v^2x_0^4 + 6Dl^2n^4v^2x_0^5 - 4Dln^4v^2x_0^6 + Dn^4v^2x_0^7 - 12Dl^3R^2n^2vx_0^2$$
$$+ 36Dl^2R^2n^2vx_0^3 - 36DlR^2n^2vx_0^4 + 12DR^2n^2vx_0^5 - 6l^4R^4k_s + 24l^3R^4k_sx_0$$
$$- 18l^2R^4k_sx_0^2 - 24Dl^3R^4 + 72Dl^2R^4x_0 - 72DlR^4x_0^2 + 36DR^4x_0^3) \quad (E.13)$$

$$f2_{f4} = \frac{-1}{2\chi}(Dl^3n^4v^2x_0^4 - 3Dl^2n^4v^2x_0^5 + 3Dln^4v^2x_0^6 - Dn^4v^2x_0^7 + 6Dl^3R^2n^2vx_0^2$$
$$- 18Dl^2R^2n^2vx_0^3 + 24DlR^2n^2vx_0^4 - 12DR^2n^2vx_0^5 + 12l^3R^4k_sx_0 - 18l^2R^4k_sx_0^2$$
$$+ 36DlR^4x_0^2 - 36DR^4x_0^3) \quad (E.14)$$

$$f3_{v1} = \frac{-1}{2\chi}(-3Dl^4n^3v^2x_0^2 + 15Dl^3n^3v^2x_0^3 - 27Dl^2n^3v^2x_0^4 + 21Dln^3v^2x_0^5 - 6Dn^3v^2x_0^6$$
$$- 6Dl^4R^2nv + 30Dl^3R^2nvx_0 - 78Dl^2R^2nvx_0^2 + 90DlR^2nvx_0^3 - 36DR^2nvx_0^4) \quad (E.15)$$

$$f3_{v4} = \frac{-1}{2\chi}(-3Dl^3n^3v^2x_0^3 + 12Dl^2n^3v^2x_0^4 - 15Dln^3v^2x_0^5 + 6Dn^3v^2x_0^6 - 6Dl^3R^2nvx_0$$
$$+ 24Dl^2R^2nvx_0^2 - 54DlR^2nvx_0^3 + 36DR^2nvx_0^4) \quad (E.16)$$

$$f3_{w1} = \frac{-1}{2\chi}(3Dl^4n^4v^2x_0^2 - 12Dl^3n^4v^2x_0^3 + 18Dl^2n^4v^2x_0^4 - 12Dln^4v^2x_0^5 + 3Dn^4v^2x_0^6$$
$$+ 6Dl^4R^2n^2v - 24Dl^3R^2n^2vx_0 + 54Dl^2R^2n^2vx_0^2 - 60DlR^2n^2vx_0^3$$
$$+ 24DR^2n^2vx_0^4 + 36l^2R^4k_sx_0 - 36lR^4k_sx_0^2 + 36Dl^2R^4 - 72DlR^4x_0$$
$$+ 36DR^4x_0^2) \quad (E.17)$$

$$f3_{w4} = \frac{-1}{2\chi}(-3Dl^2n^4v^2x_0^4 + 6Dln^4v^2x_0^5 - 3Dn^4v^2x_0^6 + 24DlR^2n^2vx_0^3 - 24DR^2n^2vx_0^4$$
$$- 36l^2R^4k_sx_0 + 36lR^4k_sx_0^2 - 36Dl^2R^4 + 72DlR^4x_0 - 36DR^4x_0^2) \quad (E.18)$$

$$f3_{f1} = \frac{-1}{2\chi}(Dl^4n^4v^2x_0^3 - 4Dl^3n^4v^2x_0^4 + 6Dl^2n^4v^2x_0^5 - 4Dln^4v^2x_0^6 + Dn^4v^2x_0^7 + 6Dl^4R^2n^2vx_0$$
$$- 24Dl^3R^2n^2vx_0^2 + 42Dl^2R^2n^2vx_0^3 - 36DlR^2n^2vx_0^4 + 12DR^2n^2vx_0^5 - 6l^4R^4k_s$$
$$+ 24l^3R^4k_sx_0 - 18l^2R^4k_sx_0^2 + 36Dl^2R^4x_0 - 72DlR^4x_0^2 + 36DR^4x_0^3) \quad (E.19)$$

$$f3_{f4} = \frac{-1}{2\chi}(Dl^3n^4v^2x_0^4 - 3Dl^2n^4v^2x_0^5 + 3Dln^4v^2x_0^6 - Dn^4v^2x_0^7 - 12Dl^2R^2n^2vx_0^3$$
$$+ 24DlR^2n^2vx_0^4 - 12DR^2n^2vx_0^5 + 12l^3R^4k_sx_0 - 18l^2R^4k_sx_0^2 + 12Dl^3R^4$$
$$- 36Dl^2R^4x_0 + 36DlR^4x_0^2 - 36DR^4x_0^3) \quad (E.20)$$

and

$$\chi = L(Dl^3n^4v^2x_0^3 - 3Dl^2n^4v^2x_0^4 + 3Dln^4v^2x_0^5 - Dn^4v^2x_0^6 + 3l^3R^4k_s + 12Dl^2R^4$$
$$- 36DLR^4x_0 + 36DR^4x_0^2) \quad (E.21)$$

**References**

1. Crisfield, M.A.: Non-linear Finite Element Analysis of Solids and Structures Volume 2: ADVANCED TOPICS. JOHN WILEY & SONS, Chichester (1997)

2. Nguyen-Xuan, H., Liu, G.R., Bordas, S., Natarajan, S., Rabczuk, T.: An adaptive singular ES-FEM for mechanics problems with singular field of arbitrary order. Comput Methods Appl Mech Eng. 253, 252–273 (2013)

3. Khoei, A.R.: Extended Finite Element Method. John Wiley & Sons, Chichester (2015)

4. Surendran, M., Natarajan, S., Bordas, S.P.A., Palani, G.S.: Linear smoothed extended finite element method. Int J Numer Methods Eng. 112(12), 1733–1749 (2017)

5. Bordas, S., Nguyen, P.V., Dunant, C., Guidoum, A., Nguyen-Dang, H.: An extended finite element library. Int J Numer Methods Eng. 71(6), 703–732 (2007)

6. Liu, G.R.: MESHFREE METHODS Moving Beyond the Finite Element Method (SECOND EDITION). Taylor & Francis Group, Boca Raton (2010)

7. Rabczuk, T., Bordas, S.P.A., Askes, H.: Meshfree Discretization Methods for Solid Mechanics. In: Encyclopedia of Aerospace Engineering. John Wiley & Sons, Chichester (2010)

8. Kienzler, R., Herrmann, G.: An Elementary Theory of Defective Beams. Acta Mech. 62, 37–46 (1986)

9. Ricci, P., Viola, E.: Stress intensity factors for cracked T-sections and dynamic behaviour of T-beams. Eng Fract Mech. 73(1), 91–111 (2006)

10. Okamura, H., Liu, H.W., Chu, C.S., Liebowitz, H.: A cracked column under compression. Eng Fract Mech.1(3), 547–564 (1969)

11. Tada, H., Paris, P.C., Irwin, G.R.: The Stress Analysis of Cracks Handbook, Third Edition. The American Society of Mechanical Engineers, New York (2000)

12. Tharp, T.M.: A finite element for edge-cracked beam columns. Int J Numer Methods Eng. 24(10), 1941–1950 (1987)

13. Yokoyama, T., Chen, M.C.: Vibration analysis of edge-cracked beams using a line-spring model. Eng Fract Mech. 59(3), 403–409 (1998)

14. Attar, M., Karrech, A., Regenauer-Lieb, K.: Free vibration analysis of a cracked shear deformable beam on a two-parameter elastic foundation using a lattice spring model. J Sound Vib. 333(11), 2359–2377 (2014)

15. Skrinar, M.: Elastic beam finite element with an arbitrary number of transverse cracks. Finite Elem Anal Des. 45(3),181–189 (2009)

16. Alijani, A., Mastan Abadi, M., Darvizeh, A., Abadi, M.K.: Theoretical approaches for bending analysis of founded Euler–Bernoulli cracked beams. Arch Appl Mech. 88(6), 875-895 (2018)

17. Mottaghian, F., Darvizeh, A., Alijani, A.: A novel finite element model for large deformation analysis of cracked beams using classical and continuum-based approaches. Arch Appl Mech. 89(2), 195–230 (2019)

18. Alijani, A., Abadi, M.K., Razzaghi, J., Jamali, A.: Numerical analysis of natural frequency and stress intensity factor in Euler–Bernoulli cracked beam. Acta Mech. 230(12), 4391–4415 (2019)

19. Alijani, A., Darvizeh, M., Darvizeh, A., Ansari, R.: Development of a semi-analytical nonlinear finite element formulation for cylindrical shells. Proc Inst Mech Eng Part C J Mech Eng Sci. 228(2), 199–217 (2014)

20. Alijani, A., Darvizeh, M., Darvizeh, A., Ansari, R.: On nonlinear thermal buckling analysis of cylindrical

shells. Thin-Walled Struct.95, 170–182 (2015)

21. Akrami, V., Erfani, S.: An analytical and numerical study on the buckling of cracked cylindrical shells. Thin-Walled Struct.119, 457–469 (2017)
22. Delale, F., Erdogan, F.: Application of the Line-Spring Model to a Cylindrical Shell Containing a Circumferential or Axial Part-Through Crack. J Appl Mech. 49, 97–102 (1982)
23. Ezzat, H., Erdogan, F.: Elastic-plastic fracture of cylindrical shells containing a part-through circumferential crack. J Press Vessel Technol Trans ASME. 104(4), 323–330 (1982)
24. Moradi, S., Tavaf, V.: Crack detection in circular cylindrical shells using differential quadrature method. Int J Press Vessel Pip. 111–112, 209–216 (2013)
25. El Naschie, M.S.: A branching solution for the local buckling of a circumferentially cracked cylindrical shell. Int J Mech Sci. 16(10), 689–697 (1974)
26. Nasirmanesh, A., Mohammadi, S.: Eigenvalue buckling analysis of cracked functionally graded cylindrical shells in the framework of the extended finite element method. Compos Struct. 159, 548–566 (2017)
27. Natarajan, S., Chakraborty, S., Ganapathi, M., Subramaniam, M.: A parametric study on the buckling of functionally graded material plates with internal discontinuities using the partition of unity method. Eur J Mech - A/Solids. 44, 136–147 (2013)
28. Venkatachari, A., Natarajan, S., Ganapathi, M., Haboussi, M.: Mechanical buckling of curvilinear fibre composite laminate with material discontinuities and environmental effects. Compos Struct. 131,790–798 (2015)
29. Rajagopalan, K.: Finite element buckling analysis of stiffened cylindrical shells. A.A. Balkema, Rotterdam (1993)
30. El-Kaabazi, N., Kennedy, D.: Calculation of natural frequencies and vibration modes of variable thickness cylindrical shells using the Wittrick-Williams algorithm. Comput Struct. 104–105, 4–12 (2012)
31. Li, Y., Zhang, Y., Kennedy, D.: Random vibration analysis of axially compressed cylindrical shells under turbulent boundary layer in a symplectic system. J Sound Vib. 406, 161–180 (2017)
32. Venstel, E., Krauthammer, T.: Thin Plates and Shells; Theory, Analysis, and Application. Marcel Dekker, New York, (2001)
33. Wood, R.D., Schrefler, B.: Geometrically non-linear analysis—A correlation of finite element notations. Int J Numer Methods Eng. 12(4), 635–642 (1978)
34. Leissa, A.W., Qatu, M.S.: Vibrations of Continuous Systems. McGraw-Hill, US (2011)
35. Moazzez, K., Saeidi Googarchin, H., Sharifi, S.M.H.: Natural frequency analysis of a cylindrical shell containing a variably oriented surface crack utilizing Line-Spring model. Thin-Walled Structures. 125, 63-75 (2018)
36. Bazant, Z. P., Cedolin, L.: Stability of structures: elastic, inelastic, fracture, and damage theories. Oxford University, N.Y (1991)
37. Singer, J., Arbocz, J., Weller, T.: Buckling Experiments: Experimental Methods in Buckling of Thin-Walled Structures: Basic Concepts, Columns, Beams and Plates- Volume 1. John Wiley & Sons, New York (1998)

38.     Wriggers, P. Nonlinear Finite Element Methods. Springer, Berlin Heidelberg (2008)